\documentclass{amsart}

\usepackage{amsmath,amsfonts,amssymb}
\usepackage{latexsym, url}
\usepackage[perpage,symbol*]{footmisc} 

\newtheorem{thm}{Theorem}[section]
\newtheorem{lem}[thm]{Lemma}

\theoremstyle{definition}

\newcommand{\Z}[1]{\mathbb{Z}/#1\mathbb{Z}}

\makeatletter
\newenvironment{remark}[1][Remark]{\par
  \normalfont \topsep6\p@\@plus6\p@\relax
  \trivlist
  \item[\hskip\labelsep
        \itshape
    #1\@addpunct{:}]\ignorespaces}
    {\endtrivlist\@endpefalse}
\makeatother

\begin{document}

\title[Rationality of the zeta function on finite abelian $p$-groups]{Rationality of the zeta function of the subgroups of
  abelian $p$-groups}

\author{Olivier Ramar\'e}

\address{%
CNRS / Institut de Math\'ematiques de Marseille\\
Aix Marseille Universit\'e, U.M.R. 7373\\
Site Sud, Campus de Luminy, Case 907\\
13288 MARSEILLE Cedex 9, France}

\email{olivier.ramare@univ-amu.fr}




\date{August 16, 2015}


\dedicatory{To Ramdinmawia Vanlalngaia}

\begin{abstract}
  Given a finite abelian $p$-group $F$, we prove an efficient recursive
  formula for $\sigma_a(F)=\sum_{\substack{H\leq F}}|H|^a$ where
  $H$ ranges over the subgroups of $F$. We infer from this formula that the
  $p$-component of the corresponding zeta-function on groups of $p$-rank
  bounded by some constant $r$ is rational with a simple denominator. We
  also provide two explicit examples in rank $r=3$ and $r=4$ as well as a
  closed formula for $\sigma_a(F)$.
\end{abstract}

\maketitle

\section{Introduction}

The subgroups of finite abelian $p$-groups have been intensively studied. An
early paper of G.~Birkhoff establishes in \cite[Theorem 8.1]{Birkhoff*34}
material to count the number of subgroups of a given type; the version
given in~\cite[(1)]{Butler*87} is surely easier to grasp. To fix the notation,
our $p$-groups will be of rank below some fixed~$r$ and are thus
isomorphic to a product
\begin{equation}
  \label{defSNF}
  F=\Z{p^{f_1}}\times\Z{p^{f_1+f_2}}\times\cdots\times\Z{p^{f_1+f_2+\cdots+f_r}}
\end{equation}
where $f_i$ are non-negative integers. We write $F=[p;
f_1,f_2,\cdots,f_r]$. The \emph{type} of $F$ is the partition
$(f_1+\cdots+f_r,f_1+\cdots+f_{r-1},\cdots, f_1)$. The \emph{type} of
a subgroup $H$ is its type as an abstract group while its
\emph{cotype} is the type of $F/H$.  In the fifties, P.~Hall
considered the numbers
$g^{\lambda}_{\mu,\nu}(p)$ 
of subgroups of type $\mu$ and cotype~$\nu$ in a group of
type~$\lambda$ and used them as multiplication constants to form the
now called \emph{the Hall algebra}.  The combinatorial aspects have
been further developed by T.~Klein in \cite{Klein*68}, in the
milestone book of I.G.~MacDonald \cite{Macdonald*95}, and by L.~Butler
in \cite{Butler*91} and \cite{Butler*94} concerning the poset formed
by these subgroups and the inclusion. This short bibliography is by no
means complete!  Two closely related fields of enumerative algebra
concern the number of subgroups of not-especially abelian groups,
e.g. Y.~Takegahara in~\cite{Takegahara*00}, and the number of
subgroups of a given index in a fixed group, see
\cite{Grunewald-Segal-Smith*88} by F.J.~Grunewald, D.~Segal \& G.C.~Smith and
the book \cite{Lubotzky-Segal*03} by A.~Lubotzky \& D.~Segal.

Given a finite abelian group $F$ and a complex number $a$, we concentrate in
this paper on the
counting function
\begin{equation}
  \label{defsigmaa}
  \sigma_a(F)=\sum_{\substack{H\leq F,\\ \text{$H$ subgroup}}}|H|^a.
\end{equation}
We obviously have
$\sigma_a(F_1F_2)=\sigma_a(F_1)\sigma_a(F_2)$ whenever $F_1$ and $F_2$ have
coprime orders, therefore reducing the study of~$\sigma_a$ to the case of
$p$-groups.

Despite the wealth of work on the question and our restriction to finite
abelian groups, it is difficult to get formulae for $\sigma_a$ that are not (very)
intricate. Still we know that, once the type of $F$ is fixed, say equal to
$\lambda$, the value $\sigma_a(F)$ is a polynomial in $p$ and $q=p^a$ since,
by using the Hall polynomials $g^{\lambda}_{\mu,\nu}$, we have
\begin{equation*}
  \sigma_a(F)=
  \sum_{\substack{\mu, \nu}}g^{\lambda}_{\mu,\nu}(p)p^{a(\mu_1+\mu_2+\mu_3+\cdots)},
\end{equation*}
the sum being over all possible choices of $\mu$ and $\nu$.  By
combining the expression given in~\cite[(1)]{Butler*87} and the
development of the $p$-binomial coefficient given in \cite[Theorem
6.1]{Kac-Cheung*02}, we even conclude that $\sigma_a(F)$ is a
polynomial in $p$ and $q$ with integer non-negative coefficients.  The
main novelty of our study is the ``simple'' recursion formula given in
Theorem~\ref{mainrec}. As an interesting consequence, the relevant
generating series is shown to be rational;
we even provide a fully explicit
formula.

\begin{thm}
  \label{expl1}
  We have
  \begin{multline*}
    \sum_{f_1,f_2,\cdots,f_r\ge0}\mkern-20mu
    \sigma_a([p; f_1,f_2,\cdots,f_r])
    X_1^{f_1}\cdots X_r^{f_r}
    =
    \\
    \sum_{(\epsilon_k)\in\{-1,1\}^r}\prod_{t=1}^{r}
    \frac{-\epsilon_{t}p^{\epsilon^*_{t}(a+\sum_{h=t+1}^{r}\epsilon_{h})}}{p^{a+\sum_{h=t+1}^{r}\epsilon_{h}}-1}
    \frac{1}{1-p^{(a+r-t+1)\sum_{h=t}^{r}\epsilon^*_{h}-(\sum_{h=t}^{r}\epsilon^*_{h})^2}X_{t}}
  \end{multline*}
  where  $\epsilon^*=(1-\epsilon)/2$. This series belongs in particular to
  $\mathbb{Q}(p, p^a, X_1,\cdots, X_r)$ and a denominator is given in~\eqref{defB}.
\end{thm}
This formula appears already in the unpublished thesis of G. Bhowmik
\cite[Section IX]{Bhowmik*97}, with whom I collaborated at that time.
By ``a denominator'', we mean a polynomial by which we can multiply
our series to fall in~ $\mathbb{Q}[p,p^a,X_1,\cdots,X_r]$. No
minimality is assumed.  The dependence in $p^a=q$ is maybe better
explained by modifying~\eqref{defsigmaa} in case of a $p$-group~into
\begin{equation}
  \label{defsigmaabis}
  \sigma_a(F)=\sum_{\substack{H\leq F,\\ \text{$H$
        subgroup}}}q^{\log|H|/\log p}.
\end{equation}
We infer from Theorem~\ref{expl1} the rather compact closed formula
\eqref{closedform} for $\sigma_a(F)$.


\subsection*{A detour in integer matrices arithmetic}

The proof below being easier to understand in the framework of integer matrices,
let us present this  hundred years old field called
 \emph{Noncommutative Number Theory} by
L.N.~Vaserstein in~\cite{Vaserstein*87}.
The
book~\cite{MacDuffee*46} of C.C.~MacDuffee contains already, in this context, a notion of gcd
and lcm
that is till under scrutiny, see~\cite{Thompson*86} by R.C.~Thompson.  A
founding result is that, when decomposing a non-singular integer matrix $M$ 
 as
a product of two integer matrices $M=AB$, the number of right-classes of $A$
under the action of $SL_r(\mathbb{Z})$ is finite; $A\cdot SL_r(\mathbb{Z})$ is
then called a left-divisor of $M$. From this fact, V.C.~Nanda in
\cite{Nanda*84} and \cite{Nanda*90} introduced a convolution product
between functions of integer matrices invariant under the action of $SL_r(\mathbb{Z})$. 
This algebra is (almost immediately) isomorphic to the Hecke algebra,
see the book \cite{Krieg*90} by A.~Krieg. 
V.C.~Nanda detailed examples among
which we find an Euler totient function, the divisor function (our $\sigma_0$), and a M\"obius function.
The initial interest for this
arithmetic comes from  modular forms.

\subsection*{Back to finite abelian groups}

Any finite abelian group of rank~$r$ can be represented as a quotient
$\mathbb{Z}^r/M(\mathbb{Z}^r)$ for some non-singular integer matrix~$M$. This
correspondence is shown in~\cite{Bhowmik-Ramare*98-2} to carry through to the
subgroups that in return appear as left divisors of~$M$.
  The left-divisibility of divisors translates as the inclusion of subgroups, and the
right-complementary divisor of any left divisor $H$ of $F$ is
associated to the quotient $F/H$. In this manner, the arithmetic of
subgroups of finite abelian groups and the one of integer matrices
locally (i.e. once a home group $F$ is chosen) coincide; for instance
the Moebius function defined on the lattice of subgroups is identical to
the one defined on matrices as the convolution inverse  of the $1\!\!1$-function. More fundamentally, the
Hall algebra, the Hecke algebra and the algebra of \emph{arithmetical}
functions on integer matrices coincide. Other connections exist;
For instance, the paper \cite{Thompson*89} of R.C.~Thompson
converts T.~Klein's combinatorial result \cite{Klein*68} in terms of divisibility of invariant
factors.

\subsection*{Average results}
 
Here, the vocabularies of groups and of matrices get mixed.  As shown by
G.~Bhowmik in~\cite{Bhowmik*91}, the function $\sigma_a(F)$ taken on
average under \emph{the determinant condition} $|F|\le x$ exhibits some regularity:
when translated in terms of abelian groups, the question is to decide of the
asymptotic behavior, when $x$ goes to infinity, of
\begin{equation*}
  \sum_{|F|\le x} \sigma_a(F) /  \sum_{|F|\le x} 1,
\end{equation*}
where $F$ ranges over the finite abelian groups of rank below some
fixed $r$.  The sum $\sum_{|F|\le x} 1$ has been the subject of
numerous publications, e.g.~\cite{Robert-Sargos*06} by A.~Ivi\`c~\cite{Ivic*92},
O.~Robert \& P.~Sargos or~\cite{Liu*10} by
H.-Q.~Liu.  The average order of $\sigma_a$
is closely related to the behavior of the rather mysterious Dirichlet series
\begin{equation*}
  D_{r,a}(s)=\prod_{p\ge2}\sum_{f_1,\cdots,f_r\ge0}
  \frac{\sigma_a([p; f_1,\cdots,f_r])}{p^{(rf_1+(r-1)f_2+\cdots+f_r)s}},
\end{equation*}
the product being taken over the primes~$p$.  Its abscissa of
convergence has been determined in \cite{Bhowmik-Ramare*94} while G.~Bhowmik \& J. Wu in \cite{Bhowmik-Wu*01} exhibit a
representation of $D_{r,a}(s)$ that yields a larger domain of meromorphic
continuation. Since the $p$-factor of this
series is the case $X_t=1/p^{s(r-t+1)}$ of Theorem~\ref{expl1}, we now have a completely explicit expression.  This
series is an analog in the finite group case of the
zeta-function introduced and studied by F.J.~Grunewald, D.~Segal \&
G.C.~Smith in \cite{Grunewald-Segal-Smith*88}, though these authors
work with a fixed group and investigate the generating function
associated to the number of subgroups of a given index, as this index
varies. In our case, the subgroups are less precisely determined (we
do not fix the index) but the sum runs over a family of~groups.  We
further note that it (as well as the more general version considered
in Theorem~\ref{expl1}) has also been investigated by
V.M.~Petrogradsky in \cite{Petrogradsky*07}.

As a side-note, we mention another kind of mean-regularity that has been obtained  in~\cite{Bhowmik-Ramare*98-1}: we have
 $\sigma_0(F)=(\log|F|)^{(1+o(1))\log 2}$ for all but $o(x)$
abelian groups of order not more than~$x$. On restricting
the set to groups of rank~$r$ exactly (there are about $x^{1/r}$ such groups), we show that $\sigma_0(F)=|F|^{[r^2/4]/r}(\log
|F|)^{\xi_r+o(1)}$ for all but $o(x^{1/r})$ exceptions, where $\xi_r=(1+(-1)^r)/2$.

\bigskip

In Section~\ref{cons}, we use our method to derive two new explicit
formulae: one when $r=3$, under the determinant condition and a
general~$a$, and one when $r=4$, still under the determinant condition though
this time restricted to the case $a=0$ to keep the expression within a
reasonable size. Finally, in Section~6, we use Theorem~\ref{expl1} to derive a closed formula for $\sigma_a(F)$.

\section{Duality}

The function $\sigma_a(F)$ is defined 
in~\eqref{defsigmaa} and we propose now another expression that is surely not novel but for which there lacks an easy reference. We present a proof for the sake of completeness.
\begin{lem}
  \label{duality}
  When $F$ is a finite abelian group, we have
  \begin{equation*}
    \sigma_a(F)=\sum_{\substack{H\leq F,\\ \text{$H$ subgroup}}}|F/H|^a.
  \end{equation*}
\end{lem}
In terms of divisors of matrices, as explained in the introduction, the expression~\eqref{defsigmaa} can be seen
as summing over left-divisors while the above can be seen as summing over
right-divisors. We present an independent
proof.
\begin{proof}
  The character group $\hat{F}$ of $F$ being isomorphic to $F$, we have
  $\sigma_a(F) = \sigma_a(\hat{F})$. The following function is known to be one-to-one, see \cite[Theorem 13.2.3]{Hall*59}:
  \begin{equation*}
    V_{F\rightarrow \hat{F}}:\begin{array}[t]{rcl}
      \text{subgroups of $F$}&\rightarrow&\text{subgroups of $\hat{F}$}\\
      H&\mapsto& H^\perp=\{\chi/\chi_{|_H}=1\}.
    \end{array}
  \end{equation*}
  It is further
  classical that $H^\perp\cong F/H$.  As a consequence, we find that
  \begin{equation*}
    \sigma_a(F)=\sum_{\substack{H\leq G,\\ \text{$H$
          subgroup}}}|H^\perp|^a
    =\sum_{\substack{H\leq F,\\ \text{$H$ subgroup of $F$}}}|F/H|^a
  \end{equation*}
  as wanted.
\end{proof}

\section{Recursion formulae}

This section is the heart of the whole paper. The next theorem together with
Lemma~\ref{duality} are the only places where we input information on our
function. Once this formula is established, the remainder of the proof of
Theorem~\ref{expl1} is maybe not immediate but is essentially a matter of bookkeeping.

\begin{thm}
  \label{mainrec}
  Let $F_r$ be a finite abelian $p$-group of rank $r\ge1$ and exponent~$p^\ell$.
  Let $e_\ell$ be an element of order $p^\ell$ and let $F_{r-1}$ be a
  subgroup such that 
  $F_r=F_{r-1}\oplus\mathbb{Z}e_\ell$.
  We have
  \begin{equation*}
    (p^a-1)\sigma_a(F_r)=
    p^{a\ell+a}\,|F_{r-1}|\,\sigma_{a-1}(F_{r-1})
    -\sigma_{a+1}(F_{r-1}).
  \end{equation*}
\end{thm}

\begin{proof}
  We consider $G_r=F_{r-1}\oplus\mathbb{Z}pe_\ell$. We first prove the
  following two recursion formulae:
  \begin{equation}
    \label{rec1}
    \sigma_a(F_r)=p^a\sigma_a(G_r)+\sigma_{a+1}(F_{r-1})
  \end{equation}
  and
  \begin{equation}
    \label{rec2}
    \sigma_a(F_r)=\sigma_a(G_r)+p^{a\ell}\,|F_{r-1}|\,\sigma_{a-1}(F_{r-1}).
  \end{equation}
  A linear combination of both gives the recursion announced in the lemma.
  The first formula will come from the expression of Lemma~\ref{duality}
  \begin{equation*}
    \sigma_a(F_r)=\sum_{\substack{H\leq F_r}} |F_r/H|^a
  \end{equation*}
  while the second one will come from the initial expression
  \begin{equation*}
    \sigma_a(F_r)=\sum_{\substack{H\leq F_r}} |H|^a.
  \end{equation*}
  To do so, we split both summations according to whether
  $H$ is a subgroup of $G_r$ or not. The first case is readily handled via the
  two formulae 
  \begin{equation}
    \label{oh1}
    \sum_{\substack{H\leq F_r,\\ H\leq G_r}} |F_r/H|^a
    =p^a\sum_{\substack{ H\leq G_r}} |G_r/H|^a=p^a\sigma_a(G_r)
  \end{equation}
  and
  \begin{equation}
    \label{oh2}
    \sum_{\substack{H\leq F_r,\\ H\leq G_r}} |H|^a
    =\sum_{\substack{ H\leq G_r}} |H|^a=\sigma_a(G_r).
  \end{equation}
  The second case requires some more analysis. Let $K$ be a subgroup of
  $F_{r-1}$. We consider
  \begin{equation*}
    \Psi:
    \begin{array}[t]{rcl}
      \bigl\{H|H\not\leq G_r, H\cap F_{r-1}=K\bigr\}
      &\rightarrow&F_{r-1}/K\\
      H&\mapsto& y\mod K, \text{where\ }y\in (H-e_\ell)\cap F_{r-1}.
    \end{array}
  \end{equation*}
This function is well-defined. Indeed, the set $(H-e_\ell)\cap F_{r-1}$ is
non-empty since $H\not\leq G_r$ and thus there exists $x=f+n e_\ell\in H$
where $f\in F_{r-1}$ and $n$ is prime to~$p^\ell$. On multiplying by the inverse of
$n$ modulo~$p$, we recover an element of the form $y+e_\ell$ as
wanted. Furthermore, the class of $y$ modulo $K$ does not depend on the choice
of~$y$. For, if $y'$ also belongs to $(H-e_\ell)\cap F_{r-1}$, then
$y-y'=(y+e_\ell)-(y'+e_\ell)$ belongs to $H\cap F_{r-1}=K$. We note that
$H=\langle K,y+e_\ell\rangle$ and that this defines the reverse function
to~$\Psi$, proving that $\Psi$ is one-to-one and onto. Note that
$F_r/H\cong F_{r-1}/K$.
As a corollary, we get
\begin{equation}
  \label{oh3}
  \sum_{\substack{H\leq F_r,\\ H\not\leq G_r}} |F_r/H|^a
  = \sum_{\substack{K\leq F_{r-1}}} |F_{r-1}/K||F_{r-1}/K|^a
  =\sigma_{a+1}(F_{r-1})
\end{equation}
and
\begin{equation}
  \label{oh4}
  \sum_{\substack{H\leq F_r,\\ H\not\leq G_r}} |H|^a
  = \sum_{\substack{K\leq F_{r-1}}} p^{a\ell}|F_{r-1}/K||K|^a
  =p^{a\ell}|F_{r-1}|\sigma_{a-1}(F_{r-1}).
\end{equation}
Combining \eqref{oh1} together with \eqref{oh3} gives \eqref{rec1} while
combining \eqref{oh2} together with \eqref{oh4} gives \eqref{rec2}.
\end{proof}

\begin{remark}[Remark 1]
  The recursion formula we prove in \eqref{rec1} is already contained
in~\cite{Bhowmik*91} where a proof in terms of matrices is given. The proof
below uses the group-theoretical context, offering the advantage that we can
re-use the same scheme of proof on the dual group, giving rise
to~\eqref{rec2}. The comparison of both yields the theorem. The reader should
notice that this formula offers a very fast manner to compute $\sigma_a(F_r)$:
the recursion~\eqref{rec1} yields an algorithm of complexity
$2^{f_1+f_2+\cdots+f_r}$ while the above reduces this complexity to $2^r$.
\end{remark}
\begin{remark}[Remark 2]
  The part of the proof that involves $\Psi$ is in fact similar to \cite[Lemma
  1.3.1 (i)]{Lubotzky-Segal*03} where complements of a given subgroup are
  being counted.
\end{remark}

In case the $r=1$, Theorem~\ref{mainrec} recovers, when $a\neq0$,
the classical formula for the sum $\sigma_a(m)=\sum_{d|m}d^a$ of the $a$-th power of
the divisors of the integer~$m$:
\begin{equation*}
  \sigma_a(p^{f_1})=\frac{p^{a(f_1+1)}-1}{p^a-1}
  =\frac{q^{f_1+1}-1}{q-1}
\end{equation*}
and by continuity, $\sigma_0(p^{f_1})=f_1+1$.
We can also use an algebraic argument: the expression for $\sigma_a$ is a
polynomial in $q=p^a$ which we evaluate at $q=1$.

By the classification of finite abelian groups, 
$F_{r-1}=[p; f_1,f_2,\cdots,f_{r-1}]$. The situation is however not
so simple concerning the subgroup $G_r$ introduced in the above proof.
Indeed, when $f_r\ge1$, we
have $G_r=[p; f_1,f_2,\cdots,f_r-1]$ but this is not the case
anymore when $f_r=0$. This fact  explains the difficulties met in~\cite{Bhowmik-Ramare*94} and
in~\cite{Bhowmik-Wu*01}. The novelty of Theorem~\ref{mainrec} is that it
produces a recurrence formula that preserves this representation.

In the case $r=2$, Theorem~\ref{mainrec} gives, when $a\neq0$,
\begin{equation}
  \label{closedformr=2}
  (p^a-1)\sigma_a([p; f_1,f_2])=
  p^{a(f_1+f_2+1)+f_1}\frac{p^{(a-1)(f_1+1)}-1}{p^{a-1}-1}
  -\frac{p^{(a+1)(f_1+1)}-1}{p^{a+1}-1}.
\end{equation}
This formula is generalized in~\eqref{closedform}.

\section{Proof of Theorem~\ref{expl1}}

Let us use our recursion  to derive an explicit formula for
\begin{equation}
  \label{defQra}
  Q_{r,a}(X_1,\cdots,X_r)=\sum_{f_1,f_2,\cdots,f_r\ge0}
  \sigma_a([p; f_1,\cdots,f_r])X_1^{f_1}\cdots X_r^{f_r}
\end{equation}
where $r\ge0$ and the parameter $p$ is fixed. The exponent $p^\ell$ of
the group $\langle f_1,\cdots,
f_r\rangle$ is $p^{f_1+\cdots+f_r}$ and we recall that $F_{r-1}=\langle
f_1,f_2,\cdots,f_{r-1}\rangle$. An immediate consequence of Theorem~\ref{mainrec} reads
\begin{align}
  \label{rel}
  (p^a-1)Q_{r,a}(&X_1,\cdots,X_r)
  =
  \\ 
  &p^a Q_{r-1,a-1}(p^{a+r-1}X_1,p^{a+r-2}X_2,\cdots,p^{a+1}X_{r-1})
  \sum_{f_r\ge0}(p^aX_r)^{f_r}
  \\&-Q_{r-1,a+1}(X_1,\cdots,X_{r-1})\sum_{f_r\ge0}X_r^{f_r}.
\end{align}
We note for future reference that
\begin{equation}
  \label{initial}
  Q_{1,a}(X_1)=\frac{1}{(1-X_1)(1-p^aX_1)},\quad
  Q_{0,a}=1.
\end{equation}
The value at $r=0$ follows from the definition~\eqref{defQra}. We also check
directly that the relation~\eqref{rel} holds true also
when $r=1$, though we will not use it.
We rewrite the above, when $a\neq0$ and $r\ge1$, in the form
\begin{align*}
  Q_{r,a}(X_1,&\cdots,X_r)
  =
  \\ 
  &\frac{p^a}{(p^a-1)(1-p^aX_r)} Q_{r-1,a-1}(p^{a+r-1}X_1,p^{a+r-2}X_2,\cdots,p^{a+1}X_{r-1})
  \\&-
  \frac{1}{(p^a-1)(1-X_r)}
  Q_{r-1,a+1}(X_1,\cdots,X_{r-1}).
\end{align*}
We can reiterate this process to obtain a rational fraction, provided that the
parameter $a$ that appears does not vanish which we assume. We will argue by
continuity later. Each time we use the above formula, we change the parameter
$r$ to $r-1$, the parameter $a$ to $a+\epsilon$ where $\epsilon=\pm1$ and the
parameters $X_i$ to $p^{\epsilon^*(a+r-i)}X_i$ where
$\epsilon^*=(1-\epsilon)/2$. We furthermore multiply $Q_{r-1,a+\epsilon}$ by
$w_\epsilon(p^a,X_r)$ where
\begin{equation}
  \label{defw}
  w_\epsilon(q,Y)=
  -\epsilon\frac{q^{\epsilon^*}}{q-1}\frac{1}{1-q^{\epsilon^*}Y}.
\end{equation}
With these notations, the above relation reads
\begin{equation}
  \label{step1}
  Q_{r,a}((X_i))
  =
  \sum_{\epsilon\,\in\{-1,1\}}w_\epsilon(p^a,X_r)
  Q_{r-1,a+\epsilon}\bigl((p^{\epsilon^*(a+r-i)}X_i)_i\bigr).
\end{equation}
In this form, it is easily iterated and yields the next lemma.

\begin{lem}
  \label{explicite}
  When $r\ge1$, we have
\begin{equation*}
  Q_{r,a}((X_i))
  =
  \sum_{(\epsilon_k)}\prod_{s=0}^{r-1}
  w_{\epsilon_{r-s}}\Bigl(
  p^{a+\sum_{k=0}^{s-1}\epsilon_{r-k}},
  p^{\sum_{k=0}^{s-1}\epsilon^*_{r-k}(a+s-k+\sum_{\ell=0}^{k-1}\epsilon_{r-\ell})}X_{r-s}
  \Bigr)
\end{equation*}
where the sum runs over $(\epsilon_k)_{1\le k\le r}\in\{-1,1\}^r$
and $\epsilon^*=(1-\epsilon)/2$.
\end{lem}

\begin{proof}
  To prove the above formula, we use recursion, starting from $r=1$ where it
  is readily checked.  We employ~\eqref{step1} to get
  \begin{multline*}
    Q_{r,a}((X_i))
    =\\
    \sum_{\epsilon_r\in\{\pm1\}}w_{\epsilon_r}(p^a,X_r)\mkern-20mu
    \sum_{(\epsilon_k)_{1\le k\le r-1}\in\{\pm1\}^{r-1}}\prod_{s=0}^{r-2}
    w_{\epsilon_{r-1-s}}\Bigl(
    p^{a+\epsilon_r+\sum_{k=0}^{s-1}\epsilon_{r-1-k}},\qquad
    \\
    p^{\sum_{k=0}^{s-1}\epsilon^*_{r-1-k}(a+\epsilon_r+s-k
        +\sum_{\ell=0}^{k-1}\epsilon_{r-1-\ell})}
      p^{\epsilon_r^*(a+\epsilon_r+r-(r-1-s))}X_{r-1-s}
      \Bigr).
  \end{multline*}
  With $s'=s+1$, $k'=k+1$ and $\ell'=\ell+1$, the right-hand side reads:
  \begin{multline*}
    \sum_{\epsilon_r\in\{\pm1\}}w_{\epsilon_r}(p^a,X_r)
    \mkern-20mu\sum_{(\epsilon_k)_{1\le k\le r-1}\in\{\pm1\}^{r-1}}\prod_{s'=1}^{r-1}
    w_{\epsilon_{r-s'}}\Bigl(
    p^{a+\epsilon_r+\sum_{k'=1}^{s'-1}\epsilon_{r-k'}},
    \\
    p^{\sum_{k'=1}^{s'-1}\epsilon^*_{r-k'}(a+\epsilon_r+s'-k'
        +\sum_{\ell'=1}^{k'-1}\epsilon_{r-\ell})}
      p^{\epsilon_r^*(a+\epsilon_r+s')}X_{r-s'}
      \Bigr).
  \end{multline*}
  We transform the above expression with:
  \begin{align*}
    &{a+\epsilon_r+\sum_{k'=1}^{s'-1}\epsilon_{r-k'}}
    =
    {a+\sum_{k'=0}^{s'-1}\epsilon_{r-k'}},
    \\&
    {\sum_{k'=1}^{s'-1}\epsilon^*_{r-k'}\Bigl(a+\epsilon_r+s'-k'
      +\sum_{\ell'=1}^{k'-1}\epsilon_{r-\ell}\Bigr)}
    =
    {\sum_{k'=1}^{s'-1}\epsilon^*_{r-k'}\Bigl(a+s'-k'
      +\sum_{\ell'=0}^{k'-1}\epsilon_{r-\ell}\Bigr)},
    \\&
    {\sum_{k'=1}^{s'-1}\epsilon^*_{r-k'}\Bigl(a+s'-k'
      +\sum_{\ell'=0}^{k'-1}\epsilon_{r-\ell}\Bigr)}
    \ +\ {\epsilon_r^*(a+\epsilon_r+s')}
    \\&\hspace{140pt}=
    \sum_{k'=0}^{s'-1}\epsilon^*_{r-k'}\Bigl(a+s'-k'
      +\sum_{\ell'=0}^{k'-1}\epsilon_{r-\ell}\Bigr).
  \end{align*}
  The factor $w_{\epsilon_r}(p^a,X_r)$ gets readily incorporated in the
  product over $s$ from $1$ to $r-1$ as the value for $s=0$. This completes
  the proof.
\end{proof}

On using the definition given by~\eqref{defw} on the expression given by
Lemma~\ref{explicite}, we get a fully explicit formula:
\begin{multline*}
  Q_{r,a}((X_i))
  =
  \\
  \sum_{(\epsilon_k)\in\{\pm1\}^r}\prod_{s=0}^{r-1}
 \frac{-\epsilon_{r-s}p^{\epsilon^*_{r-s}(a+\sum_{k=0}^{s-1}\epsilon_{r-k})}}{p^{a+\sum_{k=0}^{s-1}\epsilon_{r-k}}-1}
  \frac{1}{1-p^{\sum_{k=0}^{s}\epsilon^*_{r-k}(a+s-k+\sum_{\ell=0}^{k-1}\epsilon_{r-\ell})}X_{r-s}}.
\end{multline*}
Some beautification is called for. We first notice that
$-k+\sum_{\ell=0}^{k-1}\epsilon_{r-\ell}=\sum_{\ell=0}^{k-1}(\epsilon_{r-\ell}-1)=-2\sum_{\ell=0}^{k-1}\epsilon^*_{r-\ell}$,
yielding that $Q_{r,a}((X_i))$ is equal to 
\begin{multline*}
  \sum_{(\epsilon_k)\in\{\pm1\}^r}\prod_{s=0}^{r-1}
  \frac{-\epsilon_{r-s}p^{\epsilon^*_{r-s}(a+\sum_{k=0}^{s-1}\epsilon_{r-k})}}{p^{a+\sum_{k=0}^{s-1}\epsilon_{r-k}}-1}
  \frac{1}{1-p^{\sum_{k=0}^{s}\epsilon^*_{r-k}(a+s-2\sum_{\ell=0}^{k-1}\epsilon^*_{r-\ell})}X_{r-s}}.
\end{multline*}
The indices of shape $r-s$, $r-k$ and $r-\ell$ were useful for the recursion,
but introduce now a useless level of complexity. We set $t=r-s$, $h=r-k$ and
$g=r-\ell$ and get, for $Q_{r,a}((X_i))$, the expression
\begin{equation*}
  \sum_{(\epsilon_k)\in\{\pm1\}^r}\prod_{t=1}^{r}
  \frac{-\epsilon_{t}p^{\epsilon^*_{t}(a+\sum_{h=t+1}^{r}\epsilon_{h})}}{p^{a+\sum_{h=t+1}^{r}\epsilon_{h}}-1}
  \frac{1}{1-p^{\sum_{h=t}^{r}\epsilon^*_{h}(a+r-t-2\sum_{g=h+1}^{r}\epsilon^*_{g})}X_{t}}.
\end{equation*}
The proof of Theorem~\ref{expl1} is almost complete. We only need to use the
identity
\begin{equation*}
  \sum_{h=t}^{r}\epsilon^*_{h}\Bigl(a+r-t-2\sum_{g=h+1}^{r}\epsilon^*_{g}\Bigr)
  =
  (a+r-t)\sum_{h=t}^{r}\epsilon^*_{h}-\Bigl(\sum_{h=t}^{r}\epsilon^*_{h}\Bigr)^2
  +\sum_{h=t}^{r}\epsilon^*_{h}
\end{equation*}
which is valid because ${\epsilon^*_{h}}^2=\epsilon^*_{h}$.

\section{Consequences on Dirichlet series}
\label{cons}

Let us investigate a possible denominator for the series $Q_{r,a}((X_i))$ of
Theorem~\ref{expl1}.  The index $t$ being fixed, for each $(\epsilon_k)$, only
one factor has the variable $X_t$. All these factors are of the shape
$1-p^{(a+r-t+1)j-j^2}X_t$ for some $j$ in $\{0, \cdots, r-t+1\}$.
A possible denominator is thus
simply 
\begin{equation}
  \label{defB}
  B_r(p,q,X_1,\cdots,X_r)=
  \prod_{t=1}^r\prod_{j=0}^{r-t+1} \bigl(1-q^jp^{(r-t+1)j-j^2}X_t\bigr)
\end{equation}
by which we means that the product $A_r=B_r(p,q,X_1,\cdots,X_r)Q_{r,a}((X_i))$ falls a priori inside $\mathbb{Q}(p,q)[X_1,\cdots,X_r]$. However the only possible remaining poles are for $q=p^b$ for some integer~$b$ and this is not possible since, when $s=2+|b|$ and $X_t=1/p^{s(r-t+1)}$, the series $D_{r,a}(s)$ is bounded. It would be helpful to get a better description of $A_r$, and at minimum show that it is prime to $B_r$. Furthermore, its coefficients are integers and thus likely to have a combinatorial expression. We will see below that these coefficients may vary in signs.

\bigskip
When we restrict our attention to the case $q=1$ (i.e. $a=0$) and
$X_t=1/p^{s(r-t+1)}$, the denominator~\eqref{defB} becomes (Careful! We have replaced $j$
by $i$ and then used $j=r-t+1$ to be able to compare with \cite[Theorem
1]{Bhowmik-Wu*01}): 
\begin{equation*}
  \prod_{j=1}^r\prod_{i=0}^{j} \bigl(1-p^{-js+ji-i^2}\bigr).
\end{equation*}
The zeta product extracted in \cite[Theorem
1]{Bhowmik-Wu*01} corresponds to $i=j/2$ when $j$ is
even and to $i=(j\pm1)/2$ when $j$ is odd.

We checked the formula given by Theorem~\ref{expl1} in the case $r=1$
with~\eqref{initial} and in the case $r=2$ with \cite[Corollary
3]{Bhowmik-Ramare*94} that we recall:
\begin{equation}
  \label{eq:2}
  B_2(p,q,X_1,X_2)\,Q_{2,a}(X_1,X_2)=1+qX_1-q(q+1)X_1X_2.
\end{equation}
In the case $a=0$ and $r=3$, the erratum \cite[(4.17)]{Bhowmik-Ramare*94-2} gives
the proper formula that we have also checked against our expression.

We finally investigated formulae for $r=3$, $r=4$ and $r=5$. The formulae are huge in
general. We can however record two new explicit
formulae to help test conjectures. When $r=3$, we can keep $a$ arbitrary and
still have a manageable 
formula under the determinant condition:
\begin{multline}
  \label{eq:1}
  B_3(p,q,X^3,X^2,X)\,Q_{3,a}(X^3,X^2,X)=
  \\
  1+qX^2+ p(q + 1)qX^3 - (q^2 + (p + 1)q + 1)qX^4
  \\-((p + 1)(q^3+1) + (p^2 + p + 1)q(q+1))qX^5 
  \\+ (q^4 + (p + 1)q(q^2+1) + (p^2 + p + 1)q^2 + 1)qX^6 
  \\- (q^2 + q + 1)pq^3X^8 + (q^3 + (p + 1)q(q+1) + 1)pq^3X^9 
  \\+ (pq^2 + (p + 1)q + p)pq^4X^{10} 
  -(pq^3 + (p + 1)q(q+1) + p)pq^4X^{11}.
\end{multline}
We have used a GP-Pari \cite{PARI-GP} script to run the computations based
on~\eqref{step1} rather than on Theorem~\ref{expl1}. Since we know a possible denominator, we have used a data structure of the form [Numerator,
Denominator-Vector], where Denominator-Vector was a list of triplets $[u,v,k]$
meaning that the denominator was the product of $1-p^uq^vX_k$ taken over all
the triplets of the list. The addition of any two such structures is readily
handled.  The computations took essentially no time, while a brute force
algorithm using Theorem~\ref{expl1} and relying on the arithmetic of rational fractions was taking a very long time when $r=4$.
We checked that the final minimal denominator was indeed $B_r(p,q,X_1,\cdots,X_r)$.
When $r=4$, we use the determinant condition and stick to $a=0$ to get
\begin{multline}
  B_4(p,1,X^4,X^3,X^2,X)\,Q_{4,0}(X^4,X^3,X^2,X)=
  \\
  (7p^3 + 5p^2 + 8p + 4)p^6X^{26} 
  - (6p^3 + 4p^2 + 6p + 2)p^6X^{25}
  \\- (5p^3 + 10p^2 + 9p + 8)p^6X^{24} 
  + (6p^4 - 10p^3 - 4p^2 - 4p -2)p^4X^{23} 
  \\+ (5p^4 + 15p^3 + 7p^2 + 8p + 1)p^4X^{22} 
  + (4p^5 + 6p^3 +12p^2 - 2p + 4)p^4X^{21} 
  \\- 3(3p^2 + p + 1)p^4X^{20} 
  - 2(2p^6 +3p^5 + 12p^4 + 5p^3 + 5p^2 + 2p + 1)pX^{19} 
  \\+ (p^7 + 3p^6 + 12p^5+ 23p^4 + 6p^3 + 8p^2 + 3p + 1)pX^{18} 
  \\+ 2(p^8 + 2p^7 + 5p^6 + 7p^5 + 18p^4 + 18p^3 + 10p^2 + 6p + 2)pX^{17} 
  \\- (3p^8 + 6p^7 +15p^6 + 16p^5 + 16p^4 + 18p^3 - 4p^2 + p - 2)pX^{16} 
  \\- (2p^6 + 6p^5 + 16p^4 + 18p^3 + 24p^2 + 24p + 4)p^2X^{15} 
  \\+ (3p^7 + p^6 + p^5- 20p^4 - 13p^3 - 21p^2 - 26p - 9)pX^{14} 
  \\+ 2(3p^6 + 4p^5 + 5p^4 +4p^3 + 4p^2 + 4p - 3)pX^{13} 
  \\+ (9p^5 + 11p^4 + 23p^3 + 15p^2 +12p + 10)pX^{12} 
  \\- 6(p^2 + p - 1)p^3X^{11} 
  - (3p^5 + 2p^4 +8p^3 - p^2 + 5p + 9)X^{10} 
  \\+ 4(p^5 + p^3 + 2p + 2)X^9 
  + (6p^4 + 4p^3+ 9p^2 + 2p + 7)X^8 
  \\- 2(p^3 + 3p^2 + 2p + 4)p^2X^7 
  - (3p^2 +2p + 3)pX^6 
  \\- 2(2p^2 + 3p + 4)X^5 
  + (p^2 + 4p + 2)pX^4 +2pX^3 + X^2 + 1
\end{multline}
This expression shows that the polynomial in $p$ in front of each power of $X$
is not a sum of monomials of constant signs, as could have been thought
from the expressions for $r\le 3$. We computed similarly the $p$-factor
$Q_{5,0}(X^5,X^4,X^3,X^2,X)$ and obtained a quotient of a polynomial in
$\mathbb{Z}[p,X]$ of degree $50$ in $X$ and $25$ in $p$, the largest monomial
being $11p^{25}X^{50}$, by~$B_5(p,1,X^5,X^4,X^3,X^2,X)$ as expected.

\section{A closed formula}

We exploit Theorem~\ref{expl1} to express $\sigma_a([p;f_1,\cdots,f_r])$. We use the expansion
\begin{equation*}
  \frac{1}{1-q^{au}p^{v}X_t}
  =
  \sum_{f_t\ge0}q^{auf_t}p^{vf_t}X_t^{f_t}
\end{equation*}
to find that
\begin{multline*}
    \sum_{f_1,f_2,\cdots,f_r\ge0}\mkern-20mu
    \sigma_a([p; f_1,f_2,\cdots,f_r])
    X_1^{f_1}\cdots X_r^{f_r}
    =
    \\
    \sum_{(\epsilon_k)\in\{-1,1\}^r}\prod_{t=1}^{r}
    \frac{-\epsilon_{t}q^{\epsilon^*_{t}}p^{\epsilon^*_{t}\sum_{h=t+1}^{r}\epsilon_{h}}}{qp^{\sum_{h=t+1}^{r}\epsilon_{h}}-1}
    \\\prod_{t=1}^{r}\sum_{f_t\ge0}q^{\sum_{h=t}^{r}\epsilon^*_{h}f_t}p^{((r-t+1)\sum_{h=t}^{r}\epsilon^*_{h}-(\sum_{h=t}^{r}\epsilon^*_{h})^2)f_t}X_t^{f_t}.
  \end{multline*}
On identifying the coefficients, we find that
\begin{equation}
\label{closedform}
\fbox{$\displaystyle\quad\begin{array}{l}
  \displaystyle\rule[10pt]{0pt}{10pt}\sigma_a([p; f_1,f_2,\cdots,f_r]) = 
\sum_{(\epsilon_k)\in\{-1,1\}^r}
  \prod_{t=1}^{r}\frac{-\epsilon_{t}q^{\epsilon^*_{t}}p^{\epsilon^*_{t}\sum_{h=t+1}^{r}\epsilon_{h}}}{qp^{\sum_{h=t+1}^{r}\epsilon_{h}}-1}
\\\hskip2cm\rule[-10pt]{0pt}{30pt}\displaystyle q^{\sum_{t=1}^{r}\sum_{h=t}^{r}\epsilon^*_{h}f_t}p^{\sum_{t=1}^{r}((r-t+1)\sum_{h=t}^{r}\epsilon^*_{h}-(\sum_{h=t}^{r}\epsilon^*_{h})^2)f_t}.
\end{array}\quad$}
\end{equation}

As a mean of verification, we note that~\eqref{defsigmaabis} gives us $\sigma_a([p; f_1,f_2,\cdots,f_r])=1$ when $q=0$. In the expression above, the only contribution when $q=0$ occurs when $\epsilon^*_t=0$ for every $t\in\{1,\cdots,r\}$, i.e. when $\epsilon_t=1$   for every $t\in\{1,\cdots,r\}$. In this case $\sum_{h=t}^{r}\epsilon^*_{h}=r-t+1$, and the above formula gives the value~1 as required.
The case $r=2$ is a given in~\eqref{closedformr=2}.

\bibliographystyle{plain}

\end{document}